\pdfoutput=1
\RequirePackage{ifpdf}
\ifpdf 
\documentclass[pdftex]{sigma}
\else
\documentclass{sigma}
\fi

\numberwithin{equation}{section}

\newtheorem{Theorem}{Theorem}[section]
\newtheorem{Lemma}[Theorem]{Lemma}
\newtheorem{Corollary}[Theorem]{Corollary}
\newtheorem{Proposition}[Theorem]{Proposition}
{ \theoremstyle{definition}
\newtheorem{Definition}[Theorem]{Definition}
\newtheorem{Example}[Theorem]{Example}
\newtheorem{Remark}[Theorem]{Remark}}

\newcommand{\g}{\mathfrak{g}}
\newcommand{\m}{\mathfrak{m}}
\newcommand{\h}{\mathfrak{h}}

\begin{document}

\newcommand{\arXivNumber}{1404.0720}

\allowdisplaybreaks

\renewcommand{\PaperNumber}{069}

\FirstPageHeading

\ShortArticleName{Harmonic Maps into Homogeneous Spaces}

\ArticleName{Harmonic Maps into Homogeneous Spaces\\ According to a
Darboux Homogeneous Derivative}

\Author{Alexandre J.~SANTANA~$^\dag$ and Sim\~{a}o N.~STELMASTCHUK~$^\ddag$}

\AuthorNameForHeading{A.J.~Santana and S.N.~Stelmastchuk}

\Address{$^\dag$~Mathematics Department, State University of Maringa (UEM), 87020-900 Maringa, Brazil}
\EmailD{\href{mailto:ajsantana@uem.br}{ajsantana@uem.br}}

\Address{$^\ddag$~Mathematics Department, Federal University of
Parana (UFPR), 86900-000 Jandaia do Sul,\\
\hphantom{$^\ddag$}~Brazil}
\EmailD{\href{mailto:simnaos@gmail.com}{simnaos@gmail.com}}

\ArticleDates{Received February 10, 2015, in f\/inal form August 24, 2015; Published online August 28, 2015}

\Abstract{Our purpose is to use a Darboux homogenous derivative to
understand the harmonic maps with values in homogeneous space. We
present a characterization of these harmonic maps from the geometry
of homogeneous space. Furthermore, our work covers all type of
invariant geometry in homogeneous space.}


\Classification{53C43; 53C30; 58E20; 60H30}

\section{Introduction}

Harmonic maps into Lie groups and homogeneous space has received
much  attention due to importance of Lie groups and
homogeneous spaces in mathematics and other Sciences.  In this
context, a~useful tool is the Darboux derivative, which provides a~kind of linearization in the sense to transfer the problem of Lie
group to Lie algebra. Considering a smooth manifold~$M$, a~Lie group
$G$ and its Lie algebra~$\g$ then the Darboux derivative of a smooth
map $F\colon M \rightarrow G$ is a~$\g$-valued 1-form given by $\alpha_F
= F^*\omega$, where $\omega$ is the Maurer--Cartan form on~$G$. It is
well-known that from $\alpha_F$ it is possible to establish
conditions such that $F$ is harmonic.

This method has been used to obtain important results about harmonic
maps on Lie groups and homogeneous space, see for example
Dai--Shoji--Urakawa~\cite{DaiY/ShojiM/UrakawaH}, Higaki~\cite{Higaki},
Khemar \cite{Khemar}, Dorf\-meis\-ter--Inoguchi--Kobayashi~\cite{DorfmeisterJ/InoguchiJ/KobayashiS}, Sharper~\cite{Sharpe},
Uhlenbeck~\cite{Uhlenbeck}. Finally, interesting theorems and nice
examples are founded in   two papers of Urakawa (see~\cite{Urakawa2}
and~\cite{Urakawa3}).

Our purpose is to study the harmonic maps using the direct relation between the reductive homogeneous space~$G/H$ and ${\mathfrak m}$, where
${\mathfrak m}$ is the horizontal component of the Lie algebra ${\mathfrak g}$ of~$G$. Currently, the main technique is based on the lift
of the smooth map $F\colon M \rightarrow G/H$ to a~smooth map of~$M$ into $G$.

About our main result, let $G$ be a Lie group, $H\subset G$ a closed Lie subgroup of $G$, $\mathfrak{g}$ and~$\mathfrak{h}$ the Lie algebras of~$G$ and~$H$, respectively. Consider~$G/H$ with a $G$-invariant connection $\nabla^{G/H}$. Our main tools in this context are the Maurer--Cartan form for
homogeneous space (see Burstall--Ranwnsley \cite{BurstallF/RanwnsleyJ}), denoted by $\mu$, and the Darboux homogeneous derivative.  The Darboux homogeneous derivative of a smooth map $F\colon M \rightarrow G/H$ is def\/ined as the $\m$-valued 1-form given by $\mu_F  =F^*\mu$. Hence it follows our main result:

\begin{theorem*}
    Let $\nabla^{G/H}$ and $\nabla^{G}$ be connections on $G/H$ and $G$, respectively. Suppose that $\pi\colon  G \rightarrow G/H$ is an affine submersion with horizontal distribution. Let $(M,g)$ be a Riemannian manifold and $F\colon M \rightarrow G/H$ a smooth map. Then $F$ is harmonic map if and only if
    \begin{gather}\label{theoremintroduction}
      d^{*}\mu_{F}  - \operatorname{tr}\big(\mu_{F} ^* \nabla^{G/H}\big) = 0.
    \end{gather}
\end{theorem*}

It is direct from equation~(\ref{theoremintroduction}) that the
harmonicity of smooth maps with values in $G/H$ depends only on
structure of $G/H$ and of the invariant geometry given by~$\nabla^{G/H}$. Furthermore, this is an analogous linearization  to
the Lie group case.

It is well known that the $G$-invariant connection $\nabla^{G/H}$ is
associated to an ${\rm Ad}(H)$-invariant bilinear form $\beta\colon \m
\times \m \rightarrow \m$. Thus assuming that $\beta$ is
skew-symmetric we conclude that a smooth map $F$ is harmonic if and
only if $d^{*}\mu_{F} = 0$. When $\beta$ has a symmetric part, and
this is context of Riemmanian homogeneous space, the study becomes
more interesting. But here we study just the special case of ${\rm
SL}(n, {\mathbb R})/{\rm SO}(n,\mathbb{R})$. In fact, we give
applications of the above theorem to the homogeneous space ${\rm
SL}(n, {\mathbb R})/{\rm SO}(n,\mathbb{R})$. First, we consider ${\rm
SL}(n, {\mathbb R})/{\rm SO}(n,\mathbb{R})$ as symmetric space. Second,
we present an application with a Fisher $\alpha$-connection in ${\rm
SL}(n, {\mathbb R})/{\rm SO}(n,\mathbb{R})$ (see also Fernandes--San
Martin~\cite{FernandesM/San-MartinL,FernandesM/San-MartinL2}). We obtain an explicit solution of
geodesics with values in the nilpotent group of the Iwasawa
decomposition of ${\rm SL}(n, {\mathbb R})/{\rm SO}(n,\mathbb{R})$ with
respect to a~special case of the $\alpha$-connections.

About the structure of this paper, in Section~\ref{section2} we give a brief summary of some concepts of connection on homogeneous spaces. Then in Section~\ref{section3} we prove our main result (see the above theorem). Finally, in Section~\ref{section4} we present some applications of the main results in case of homogeneous spaces~${\rm SL}(n, {\mathbb R})/{\rm SO}(n,\mathbb{R})$.

\section{Connections on homogeneous spaces}\label{section2}

In this section we introduce the notations and results about
homogeneous spaces that will be necessary in the sequel. We begin by
introducing the kind of homogeneous space that we will work. For
more details see Helgason~\cite{Helgason}.

Let $G$ be a Lie group and $H\subset G$ a closed Lie subgroup of $G$. Denote by $\mathfrak{g}$ and $\mathfrak{h}$ the Lie algebras of $G$ and $H$, respectively. We assume that the homogeneous space $G/H$ is reductive, that is, there is a subspace $\mathfrak{m}$ of $\mathfrak{g}$ such that $\mathfrak{g} = \mathfrak{h} \oplus \mathfrak{m}$ and ${\rm Ad}(H)(\m)\subset \m$. Take $\pi \colon G \rightarrow G/H$ the natural projection. It is well-known that there exists a neighborhood $U \subset \m$ of $0$ such that the smooth map $\phi\colon U \rightarrow \phi(U)$ def\/ined by $u \mapsto \pi(\exp u)$ is a dif\/feomorphism. Denote $N = \phi (U)$.

For each $a \in G$ we denote the left translation as $\tau_{a}\colon G/H \rightarrow G/H$ where $\tau_{a}(gH) = agH$.  Furthermore, if $L_{a}$ is the left translation on $G$, then
\begin{gather*}
  \pi \circ L_{a} = \tau_{a} \circ \pi.
\end{gather*}
As the left translation $L_{g}$ is a dif\/feomorphism we have
\begin{gather*}
  T_{g}G = (L_{g})_{*e}\mathfrak{h} \oplus (L_{g})_{*e}\mathfrak{m}, \qquad \forall\,  g \in G.
\end{gather*}
Denoting $TG_{\mathfrak{h}}:= \{ (L_{g})_{*e}\mathfrak{h}, \, \forall\, g \in G \}$ and $TG_{\mathfrak{m}}:= \{ (L_{g})_{*e}\mathfrak{m},\,\forall\, g \in G \}$ it follows that
\begin{gather}\label{decompositionTG}
  TG = TG_{\mathfrak{h}} \oplus TG_{\mathfrak{m}}.
\end{gather}
The horizontal projection of $TG$ into $TG_{\mathfrak{m}}$ is written as~$\mathbf{h}$.

From above equality we have the following facts about theory of
connection. First, we can view $\pi\colon G \rightarrow G/H$ as a~$H$-principal f\/iber bundle. Furthermore, Theorem~11.1 in~\cite{KobayashiS/NomizuK} shows that the principal f\/iber bundle
$G(G/H,H)$ has the vertical part of the Maurer--Cartan form on~$G$,
which is denoted by~$\omega$ and is a~connection form with respect
to decomposition~(\ref{decompositionTG}). It is clear that~$TG_{\m}$
is a~connection in~$G(G/H,H)$. Thus we denote the horizontal lift
from~$G/H$ to~$G$ by~$\mathcal{H}$.

Take $A \in \m$. The left invariant vector f\/ield $\tilde{A}$ on~$G$ is denoted by $\tilde{A}(g) = L_{g*}A$ and the $G$-invariant vector f\/ield $A_{*}$ on $G/H$ is def\/ined by $A_{*} =\tau_{g*}A$. It is clear that~$\tilde{A}$ is a horizontal lift vector f\/ield of $A_{*}$ and $\pi_*(\tilde{A})= A_*$. In other words, $\tilde{A}$ and~$A_{*}$ are $\pi$-related.

Now we introduce  some geometric aspects of $G$ and $G/H$. With
intention of using the theory of principal f\/iber bundles, our idea
is to choose a good connection~$\nabla^{G}$ such that it is
horizontally  projected over $\nabla^{G/H}$. In other words, $\pi\colon G
\rightarrow G/H$ will be an af\/f\/ine submersion with horizontal
distribution (see Abe--Hasegawa~\cite{AbeN/HasewagaK} for detail of
these connections or Proposition~\ref{propconn1} below).

We choose a $G$-invariant connection $\nabla^{G/H}$ on $G/H$.
Theorem 8.1 in Nomizu \cite{Nomizu} assures the existence of a
unique ${\rm Ad}(H)$-invariant bilinear form $\beta\colon  \m \times \m
\rightarrow \m$ such that
\begin{gather*}
  \big(\nabla^{G/H}_{A_*}B_*\big)_{o} = \beta(A,B), \qquad A, B \in \m.
\end{gather*}

In the following we construct a connection $\nabla^{G}$ from $\nabla^{G/H}$ such that the projection $\pi$ is an af\/f\/ine submersion with horizontal distribution. We extend $\beta$ to a bilinear form $\alpha$ from $\g \times \g$ into~$\g$ satisfying
\begin{gather}\label{alphabetaconnection}
  \mathbf{h}\alpha(A,B) = \beta(A,B)
\end{gather}
for $A,B \in \m$. Thus, there exists a left invariant connection $\nabla^{G}$ on $G$ such that
\begin{gather*}
  \big(\nabla^{G}_{\tilde{A}}\tilde{B}\big)(e) = \alpha(A,B), \qquad A,B \in \g.
\end{gather*}

Below we show that these choices make $\pi$ an af\/f\/ine submersion with horizontal distribution.

\begin{Proposition}\label{propconn1}
  Under the above assumptions, for vector fields $X$, $Y$ on $G/H$ we have
  \begin{gather*}
    \mathbf{h}\big(\nabla^{G}_{\mathcal{H}X}\mathcal{H}Y\big) = \mathcal{H}\big(\nabla^{G/H}_{X}Y\big).
  \end{gather*}
\end{Proposition}

\begin{proof}
  Take $A, B \in \m$ and the left invariant vectors f\/ields $\tilde{A}$, $\tilde{B}$ on $\exp(U)$ and the $G$-invariant vector f\/ields $A_{*}$, $B_{*}$ on~$N$. It is clear that $\tilde{A}$, $\tilde{B}$ are horizontal, $\pi_{*}(\tilde{A})= A_{*}$ and $\pi_{*}(\tilde{B}) = B_{*}$. By construction of~$\nabla^{G}$, for $g \in \exp(U)$, we have
  \begin{gather*}
    \pi_{*}\big(\nabla^{G}_{\tilde{A}}\tilde{B}\big)(g) = \pi_{*}L_{g*} \alpha(A,B) = \tau_{g*}\pi_{*}\alpha(A,B)
    = \tau_{g*}\beta(A,B) = \big(\nabla^{G/H}_{A_{*}}B_{*}\big)(\pi(g)).
  \end{gather*}
  This gives $\mathbf{h}\big(\nabla^{G}_{\tilde{A}}\tilde{B}\big) = \mathcal{H}\big(\nabla^{G/H}_{A_{*}}B_{*}\big)$. Using properties
  of connection it is easy to show that this result holds for any vector f\/ields~$X$,~$Y$ on~$N$.
  Now taking~$X$,~$Y$ on~$G/K$ we translate they to~$N$ (by invariance of~$\nabla^{G/H}$) and the result follows.
\end{proof}

In the following sections we regard $\pi$ as an af\/f\/ine submersion with horizontal distribution following the choice given by~(\ref{alphabetaconnection}) for connections~$\nabla^{G}$ and~$\nabla^{G/H}$.

\section{Harmonic maps on homogeneous spaces}\label{section3}

In this section we prove our main result and present two examples that will be important to the applications (see next section). We begin introducing the Maurer--Cartan form in the homogeneous space~$G/H$.

Denoting by $\omega$ the Maurer--Cartan form on $G$ we def\/ine the 1-form $\mu\colon T(G/H) \rightarrow \m$ as
\begin{gather*}
 \mu(\pi(g)) = \omega_{g} \circ \mathcal{H}_{\pi(g)},
\end{gather*}
where $\mathcal{H}$ is the horizontal lift from $G/H$ to $G$ given by decomposition $\g = \h \oplus \m$.  We note that~$\mu$ depends on the direct sum $\g = \h \oplus \m$. It is clear that if $A \in \m$, then $\mu(\pi(g))(A_{*}(\pi(g)))= A$.

\begin{Remark}
  Here, we note that $\mu$  can be  seen as   Maurer--Cartan form (see  \cite[Chapter~1]{BurstallF/RanwnsleyJ}). The dif\/ference is the direction of construction, while Burstall--Ranwnsley go to $\m$ for $T(G/H)$ we go in the inverse direction.
\end{Remark}

\begin{Definition}
  The homogeneous Darboux derivative of a smooth map $F\colon M \rightarrow G/H$ is the $\m$-valued $1$-form $\mu_F = F^*\mu$ on $M$.
\end{Definition}

It is  easy to show that the Maurer--Cartan form and homogeneous Darboux derivative satisfy analogous properties of the Maurer--Cartan form and Darboux derivative in Lie groups.

Before the main result we need the following technical lemmas.

\begin{Lemma}\label{le3}
  Let $(M,g)$ be a Riemannian manifold, $G/H$ a reductive homogeneous space and $F\colon M \rightarrow G/H$ a smooth map. Then, for every $1$-form $\theta$ on $\mathfrak{m}$, it follows that $d^{*}\mu_{F}^{*}\theta =  \theta d^{*}\mu_{F}$, where~$d^{*}$ is the co-differential operator on~$M$. Moreover,
  the equality $\operatorname{tr}F^*\nabla^{G/H}\mu^{*}\theta = \theta \operatorname{tr}F^*\nabla^{G/H}\mu$ holds.
\end{Lemma}

\begin{proof}
  From def\/inition of co-dif\/ferential $d^{*}$, for any orthonormal frame f\/ield $\{e_{1},\ldots, e_{n}\}$ on $M$ in a neighborhood of $p \in M$, we have
  \begin{gather*}
    d^{*}\mu_{F}^*\theta = -\sum_{i=1}^{n} (\nabla^{g}_{e_{i}}\mu_{F}^{*}\theta)(e_{i}),
  \end{gather*}
  where $\nabla^{g}$ is the Levi-Civita connection associated to the metric $g$. By def\/inition of dual connection,
  \begin{gather*}
    d^{*}\mu_{F}^{*}\theta
 =  -\sum_{i=1}^{n}\big(e_{i}\theta(\mu_{F}(e_{i})) - \theta(\mu_{F}^{*}\nabla^{g}_{e_{i}}e_{i})\big).
  \end{gather*}
  Since $\theta\colon \mathfrak{m} \rightarrow \mathbb{R}$ is a linear map and using the def\/inition of co-dif\/ferential in the last equality, we
  have
  \begin{gather*}
    d^{*}\mu_{F}^{*}\theta = \theta\left(-\sum_{i=1}^{n}\nabla^{g}_{e_{i}}(\mu_{F}(e_{i}))\right) =
    \theta(d^{*}\mu_{F}).
  \end{gather*}
  Similarly, we  show the second equality of the lemma.
\end{proof}

\begin{Lemma}\label{changeofmu}
  Let $\nabla^{G/H}$ and $\nabla^{G}$ be connections on $G/H$ and $G$, respectively, such that $\pi \colon G \rightarrow G/H$ is an affine submersion with horizontal distribution. Let~$(M,g)$ be a Riemannian manifold and $F\colon M \rightarrow G/H$ a smooth map. Then
  \begin{gather*}
    F^*\nabla^{G/H}\mu=  - \mu_{F}^*\nabla^{G/H}.
  \end{gather*}
\end{Lemma}

\begin{proof}
  First note that $F^*\nabla^{G/H}\mu$ is a tensor on $M$. Take $X, Y \in \mathcal{X}(M)$ and suppose that $F_{*}(X), F_{*}(Y) \in N$. In the case $F_{*}(X), F_{*}(Y) \notin N$, since $\mu$ is invariant by $\tau$, we consider the neighborhood $\tau(N)$. Let $A^{1}, \ldots, A^{m}$ be a basis on
  $\m$ and $A^1_*, \ldots,  A^m_*$ a frame f\/ield associated to basis on $N$. Then $F_*X = (F_*X)_{i}A^{i}_{*}$ and $F_*Y = (F_*Y)_{j}A^{j}_{*}$. Now
  \begin{gather*}
    F^*\nabla^{G/H}\mu(X, Y)
      =   \nabla^{G/H}\mu(F_*X,F_*Y)  =  \nabla^{G/H}\mu\big((F_*X)_{i}A^{i}_*,(F_*Y)_{j}A^{j}_*\big)\\
 \hphantom{F^*\nabla^{G/H}\mu(X, Y)}{} =   ((F_*X)_{i}(F_*Y)_{j})\nabla^{G/H}\mu\big(A^{i}_*,A^{j}_*\big),
  \end{gather*}
  where we use the fact that $\nabla^{G/H}\mu$ is a tensor. From def\/initions of dual connection and $\mu$ we see that
  \begin{gather*}
    \nabla^{G/H}\mu\big(A^{i}_*,A^{j}_*\big) = A^{i}_*\mu(A^{j}_*) - \mu\big(\nabla^{G/H}_{A^{i}_*}A^{j}_*\big) =
     - \omega_{g} \circ \mathcal{H} \big(\nabla^{G/H}_{A^{i}_*}A^{j}_*\big).
  \end{gather*}
  Now Proposition~\ref{propconn1} assures that
  \begin{gather*}
    \nabla^{G/H}\mu\big(A^{i}_*,A^{j}_*\big)
      =   - \omega\big( \mathbf{h}\big(\nabla^{G}_{\tilde{A}^{i}}\tilde{A}^{j}\big)\big)
      =   - \beta\big(\omega \circ \mathcal{H}\big({A}^{i}_*\big), \omega \circ \mathcal{H}\big(A^{j}_*\big)\big)
      =   - \beta\big(\mu\big({A}^{i}_*\big), \mu\big(A^{j}_*\big)\big).
  \end{gather*}
  Note that we use the left invariant property of $\nabla^{G}$ in the above second equality. Since the connection map $\beta$ is a bilinear form, it follows
  \begin{gather*}
    F^*\nabla^{G/H}\mu(X, Y)
      =   - \beta\big(\mu\big((F_*X)_{i}A^i_*\big), \mu\big((F_*Y)_{j}A^j_*\big)\big)
      =   - \beta\big(\mu(F_*X), \mu(F_*Y)\big) \\
      \hphantom{F^*\nabla^{G/H}\mu(X, Y)}{}
      =   - \nabla^{G/H}_{\mu_{F}(X)} \mu_{F}(Y) .\tag*{\qed}
  \end{gather*}
  \renewcommand{\qed}{}
\end{proof}

Now we have the main theorem. In the proof we use freely the
concepts and notations of Catuogno~\cite{Catuogno} and Emery~\cite{Emery}.

\begin{Theorem}\label{mainteo}
  Consider the same hypotheses of Lemma~{\rm \ref{changeofmu}}. Then~$F$ is harmonic map if and only if
  \begin{gather*}
    d^{*}\mu_{F}  - \operatorname{tr}\big(\mu_{F}^* \nabla^{G/H}\big) = 0.
  \end{gather*}
\end{Theorem}

\begin{proof}
  Let $F\colon M \rightarrow G/H$ be a smooth map and $B_{t}$ a Brownian motion in $M$. From the formula to convert Stratonovich's integral to the It\^o's integral (see formula (2) in \cite{Stelmastchuk}) we deduce that
  \begin{gather*}
    \int_0^t \mu d^{G/H}F(B_{s})
      =  \int_0^t \mu \delta F(B_{s}) -\dfrac{1}{2} \int_0^t \nabla^{G/H}\mu(dF(B_s),dF(B_s))\\
\hphantom{\int_0^t \mu d^{G/H}F(B_{s})}{}  =  \int_0^t (\mu_{F}) \delta B_{s} - \dfrac{1}{2}\int_0^t \nabla^{G/H}\mu(dF(B_s),dF(B_s)).
  \end{gather*}
  Thus
  \begin{gather*}
    \int_0^t  \mu d^{G/H}F(B_{s})  =  \int_0^t  (F^{*}\mu) dB_{s}  -  \dfrac{1}{2}\int_0^t
    (d^{*}\mu_{F})(B_{s})ds - \dfrac{1}{2} \int_0^t  \nabla^{G/H}\mu(dF(B_s),dF(B_s)).
  \end{gather*}
  From the last term in above right side we have
  \begin{gather*}
    \int_0^t \nabla^{G/H}\mu(dF(B_s),dF(B_s))   =   \int_0^t \nabla^{G/H}\mu(F_{*}dB_s,F_{*}dB_s)\\
\hphantom{\int_0^t \nabla^{G/H}\mu(dF(B_s),dF(B_s))}{}
 =   \int_0^t  F^{*} \nabla^{G/H}\mu(dB_s,dB_s)   =  \int_0^t  \operatorname{tr}\big(F^{*}
    \nabla^{G/H}\mu\big)(B_{s})ds .
  \end{gather*}
  Therefore
  \begin{gather}\label{teomaineq1}
    \int_0^t \mu d^{G/H}F(B_{s}) = \int_0^t \mu_{F} dB_{s} + \dfrac{1}{2}\int_0^t\big({-}d^{*}\mu_{F} - \operatorname{tr}\big(F^{*} \nabla^{G/H}\mu\big)\big)(B_{s}) ds.
  \end{gather}
  Suppose that $F$ is an harmonic map. Thus $F(B_{t})$ is a martingale in $G/H$. On the other side, $\int_0^t (F^{*}\mu) dB_{s}$ is a real local martingale. Then from equation~(\ref{teomaineq1}) and Doob--Meyer decomposition it follows that
  \begin{gather*}
    \int_0^t \big({-}d^{*}\mu_{F} - \operatorname{tr}\big(F^{*} \nabla^{G/H}\mu\big)\big)(B_{s}) ds = 0.
  \end{gather*}
  Hence, since $B$ is arbitrary,
  \begin{gather*}
    d^{*}\mu_{F}  +  \operatorname{tr}\big(F^{*} \nabla^{G/H}\mu\big) = 0.
  \end{gather*}
  Finally, Lemma \ref{changeofmu} yields
  \begin{gather}\label{teomaineq2}
    d^{*}\mu_{F}  -  \operatorname{tr}\big(\mu_{F}^* \nabla^{G/H}\big) =0.
  \end{gather}
  Conversely, suppose that (\ref{teomaineq2}) is true, and let $\theta$ be a $1$-form on $G/H$. Here, it is suf\/f\/icient to consider the 1-forms that are invariant by $\mu$, namely, 1-forms that satisfy $\theta = \theta_{o}\circ \mu$, where $\theta_o$ is viewed as $1$-form on  $\m$. Therefore
  \begin{gather*}
    \int_0^t \theta d^{G/H}F(B_{s}) = \int_0^t \theta_{o} \circ \mu  d^{G/H}F(B_{s}).
  \end{gather*}
  Analogous computations show that
  \begin{gather*}
    \int_0^t \theta d^{G/H}F(B_{s}) = \int_0^t (F^{*}\theta) dB_{s} + \dfrac{1}{2}\int_0^t \big({-}d^{*}F^{*}
    \theta_{o}\circ \mu - \operatorname{tr}\big(F^{*} \nabla^{G/H}\theta_{o} \circ \mu\big)\big)(B_{s})ds.
  \end{gather*}
  Applying Lemma~\ref{le3} in the above equality we obtain
  \begin{gather*}
    \int_0^t \theta d^{G/H}F(B_{s}) = \int_0^t (F^{*}\theta) dB_{s} + \dfrac{1}{2}\int_0^t\theta_{o}\big({-}d^{*}\mu_{F}
    - \operatorname{tr}\big(F^{*} \nabla^{G/H}\mu\big)\big)(B_{s}) ds.
  \end{gather*}
  By hypothesis,
  \begin{gather*}
    \int_0^t \theta d^{G/H}F(B_{s}) = \int_0^t (F^{*}\theta) dB_{s}.
  \end{gather*}
  Since $B_t$ is a Brownian motion, $\int_0^t (F^{*}\theta) dB_{s}$ is a real martingale. Consequently,  $\int_0^t \theta d^{G/H}F(B_{s})$ is a real martingale and, being $\theta$ arbitrary, $F(B_{t})$ is a
     $\nabla^{G/H}$-martingale. Thus, $F$ is a  harmonic map.
\end{proof}

\begin{Example}\label{exampleriemannian}
  Let $\langle \,,\,\rangle $ be a scalar product in $\m\times \m$ which is invariant by ${\rm Ad}(H)$, that is,  $\langle {\rm Ad}(h) X, {\rm Ad}(h) Y\rangle = \langle X,Y\rangle $ for $h \in H$ and $X,Y  \in \m$.  This scalar product gives a $G$-invariant metric $\langle\langle\,,\,\rangle\rangle $ on~$G/H$. The invariant connection associated to $\langle\langle\,,\,\rangle\rangle $ is given by the connection function
  \begin{gather*}
    \beta(X,Y) = \frac{1}{2}[X,Y]_{m} + U(X,Y), \qquad X, Y \in \m,
  \end{gather*}
  where $U\colon \m\times m \rightarrow \m$ is a symmetric bilinear form def\/ined by
  \begin{gather*}
    2\langle U(X,Y),Z\rangle = \langle X, [Y,Z]_{\m}\rangle  + \langle [Z,X]_{m},Y\rangle,
  \end{gather*}
  for all $X,Y,Z \in \m$ (see Nomizu~\cite{Nomizu}). Now Theorem~\ref{mainteo} shows that the necessary and suf\/f\/icient condition for
   $F\colon  M \rightarrow G/H$ to be a harmonic map is  that
  \begin{gather*}
    d^{*}\mu_{F}  - \sum_{i=1}^{n}U(\mu_{F}(e_{i}),\mu_{F}(e_{i})) = 0,
  \end{gather*}
  where $\{ e_1, \ldots , e_n\}$ is an orthonormal frame f\/ield in~$M$.
\end{Example}

\begin{Remark}
  In \cite{DaiY/ShojiM/UrakawaH}, Dai--Shoji--Urakawa presents a study about reductive Riemannian homogeneous spaces, however the authors work with lift of smooth maps to Lie group to characterize the harmonic maps with values in homogeneous space (see Theorem~2.1 in~\cite{DaiY/ShojiM/UrakawaH}).

  Assuming that~$M$ is a Riemannian surface and that $\nabla^{G/H}$ is a   canonical  connection, that is,
  $\beta(X,Y) = 0$ or $\beta(X,Y) = \frac{1}{2}[X,Y]$, Higaki in~\cite{Higaki} obtains an equivalent result of above theorem.
  Furthermore, Khemar (see Theorem~7.1.1 in~\cite{Khemar})
  also presents an equivalent result for a Riemannian surface $M$ and a family of invariant connections  given by
  $\beta_t(X,Y) = t[X,Y]$ for $0 \leq t \leq 1$.

  In \cite{DorfmeisterJ/InoguchiJ/KobayashiS}, Dorfmeister--Inoguchi--Kobayashi develop a similar work in the context
  of Lie groups. In fact, they assume that $M$ is a~Riemannian surface  and adopt a Lie group~$G$ with a left invariant
  connection instead~$G/H$ with a $G$-invariant connection.
\end{Remark}

In the sequel, we classify the harmonic maps for a large class of $G$-invariant connections. We will work with $G$-invariant connections $\nabla^{G/H}$ given by ${\rm Ad}(H)$-invariant bilinear forms $\beta$ which satisfy
\begin{gather}\label{betanull}
  \beta(X,X) = 0,\qquad \forall\, X \in \m.
\end{gather}

It is clear that every bilinear form $\beta$ can be written as a sum
of a symmetric part and a skew-symmetric part. The condition~(\ref{betanull}) is not so restrictive, in fact we can f\/ind these
bilinear forms in the   reductive homogeneous spaces~$G/H$ that
admit a $G$-invariant metric with $U \equiv 0$ and in the Riemannian
symmetric spaces. Other examples, following Nomizu~\cite{Nomizu},
are the canonical invariant connection of the f\/irst kind,  which has
$\beta(X,Y)= \frac{1}{2}[X,Y]_{m}$, for all $X,Y \in \m$, and the
canonical invariant connection of the second kind, which has
$\beta(X,Y) = 0$, for all $X,Y \in \m$. Furthermore, every
skew-symmetric bilinear form satisf\/ies the condition~(\ref{betanull}). Summarizing,
 the next result characterizes every harmonic map under this conditions.

\begin{Corollary}\label{corollarybetanull}
  Consider as above the connections $\nabla^{G/H}$, $\nabla^{G}$ and the submersion $\pi$. Let~$M$ be a Riemannian manifold and take $G/H$ a reductive homogeneous space. If~$\beta$, associated to~$\nabla^{G/H}$, satisfies $\beta(X,X)= 0$ for all $X \in \m$, then a~map   $F\colon M \rightarrow G/H$ is harmonic if and only if
  \begin{gather*}
    d^{*}\mu_{F} = 0.
  \end{gather*}
\end{Corollary}

A direct consequence of the above corollary is that the harmonic maps depend only on~$\mu$, the geometry of $\beta$ has no inf\/luence.

Now for the next corollary, let $\{A_1, \ldots, A_r\}$ be a basis of~$\m$. In this basis we have that $\mu = \sum\limits_{k=1}^{r}\mu^{k} A_k$.
Then

\begin{Corollary}
  Under assumption of Corollary~{\rm \ref{corollarybetanull}}, a map $F\colon M \rightarrow G/H$ is harmonic if and only if
  \begin{gather*}
    d^{*}\mu_{F}^k = 0,\qquad k=0,\ldots n.
  \end{gather*}
\end{Corollary}

\begin{proof}
  Let $\{A_1, \ldots, A_r\}$ be a basis of $\m$  and $\{e_1, \ldots, e_n\}$ an orthonormal frame f\/ield on~$M$. From Lemma~\ref{le3} we see that
  \begin{gather*}
    d^{*}\mu_{F} = -\sum_{i=1}^{n}\big(e_{i}(\mu_{F}(e_{i})) -
    \mu_{F}^{*}\nabla^{g}_{e_{i}}e_{i}\big).
  \end{gather*}
  In the basis of $\m$ we have that $\mu = \sum\limits_{k=1}^{r}\mu^{k} A_k$, so
  $\mu_{F}(e_{i})= \sum\limits_{k=1}^{r} \mu_{F}^k(e_i) A_{k}$. Then
  \begin{gather*}
    d^{*}\mu_{F}
      =   -\sum_{i=1}^{n}e_{i}\left(\sum_{k=1}^{r}\mu_{F}^k (e_i) A_{k}\right) =
    \sum_{k=1}^{r}\left(-\sum_{i=1}^{n}e_{i}\big(\mu_{F}^k (e_i)\big)\right) A_{k}
      =   \sum_{k=1}^{r}\big(d^*\mu_{F}^k\big) A_{k}.\tag*{\qed}
  \end{gather*}
  \renewcommand{\qed}{}
\end{proof}

In the sequel we establish two examples that are necessary in our applications.

\begin{Example}\label{iwasa-decomp}
  Let $G$ be connected non-compact semisimple Lie group with f\/inite center and denote by $\g$ its Lie algebra. Take a Cartan decomposition $\g = \h  \oplus \mathfrak{s}$. Choose a maximal abelian subspace $\mathfrak{a} \subset \mathfrak{s}$  and denote by ${\Pi}^{+}$ the corresponding set of positive roots and $\Sigma$ the set of simple roots. Put
  \begin{gather*}
    \mathfrak{n}^{+}=\sum_{\alpha \in {\prod}^{+}} {\g}_{\alpha} ,
  \end{gather*}
  where ${\g}_{\alpha}$ stands for the $\alpha$-root space. Then we have the Iwasawa decomposition on the Lie algebra, $\g = \h \oplus \mathfrak{a} \oplus \mathfrak{n}^{+}$. We denote by $H=\exp \mathfrak{h}$, $N^{+}=\exp \mathfrak{n}^{+}$ and $A=\exp \mathfrak{a}$ the connected subgroups with corresponding Lie algebras. Therefore, it follows the Iwasawa decomposition of~$G$, $G=H \times A \times N^{+}$. Thus for a smooth map $F\colon  M \rightarrow G/H$ we have that
  \begin{gather*}
    F(x) = \big(F^1(x),F^2(x)\big) \in  A\times N^{+},\qquad \forall\, x \in M.
  \end{gather*}
  Taking a vector $ v \in T_xM$ we have that $ \mu F_{x*}(v)$ is given by  $\mu_{F}(x)(v)= \mu_{F^1}(x)(v)
   + \mu_{F^2}(x)(v) \in \mathfrak{a} \oplus \mathfrak{n}^{+}$. Then
  \begin{gather*}
    d^*\mu_{F}(x)(v)  = d^*\mu_{F^1}(x)(v) + d^*\mu_{F^2}(x)(v).
  \end{gather*}
  Under condition that $\beta(X,X)=0$, from Corollary \ref{corollarybetanull} we have that $F$ is an harmonic map if and only if
  \begin{gather*}
    d^*\mu_{F^1}(x) = 0\qquad  \textrm{and}\qquad d^*\mu_{F^2}(x)=0.
  \end{gather*}
\end{Example}

\begin{Example}\label{superficieharmonic}
  Let $M^n$ be a smooth manifold  and $G/H$ a reductive homogeneous space endowed with a connection $\nabla^{G/H}$ such that its associated  bilinear form $\beta$ satisf\/ies $\beta(X,X) = 0$. Let $F\colon  M^n \rightarrow G/H$ be a smooth map and $\phi\colon \Omega \rightarrow M^n$ a local coordinate system.
   For simplicity, we write $F\colon\Omega \rightarrow G/H$, where $\Omega$ is an open set in~${\mathbb R}^n$. Then
  \begin{gather*}
    d^*\mu_{F} = \sum_{i=1}^n \frac{\partial}{\partial x^i}\mu\left(\frac{\partial F}{\partial x^i} \right).
  \end{gather*}
  Hence, from Corollary \ref{corollarybetanull} we have that   $F$ is an harmonic map if and only if
  \begin{gather*}
    \sum_{i=1}^n \frac{\partial}{\partial x^i}\mu\left(\frac{\partial F}{\partial x^i} \right) = 0.
  \end{gather*}
\end{Example}

\section{Applications}\label{section4}

In this section we present a study of harmonic maps with values in
the homogeneous spaces  ${\rm SL}(n,\mathbb{R})/{\rm
SO}(n,\mathbb{R})$. First we treat the symmetric space  ${\rm
SL}(n,\mathbb{R})/{\rm SO}(n,\mathbb{R})$ with canonical Riemannian
connection $\beta(X,Y) = \frac{1}{2}[X,Y]_m$. In the second part we
consider the  connection given by $\beta(X,Y) =
\left(\frac{XY+YX}{2} - \frac{\operatorname{tr}(XY)}{n}\operatorname{Id}\right)$, that is a
special case of Fisher $\alpha$-connections studied, for example, in~\cite{FernandesM/San-MartinL} and~\cite{FernandesM/San-MartinL2}.

Using the Iwasawa decomposition of ${\rm SL}(n,\mathbb{R})$ we have an explicit description of its homogeneous spaces. Hence, as a particular case of the Example~\ref{iwasa-decomp}, we consider the following construction for the Iwasawa decomposition of ${\mathfrak{sl}}(n, {\mathbb R})$, the Lie algebra of ${\rm SL}(n,\mathbb{R})$. Take the Cartan decomposition given by ${\mathfrak{sl}}(n, {\mathbb R})={\mathfrak{so}}(n,{\mathbb R}) \oplus {\mathfrak s}$, where
\begin{gather*}
  {\mathfrak{so}}(n,{\mathbb R})= \big\{ X \in {\mathfrak{sl}}(n, {\mathbb R})\colon  X^{t}=-X \big\} \qquad
  \mbox{and} \qquad {\mathfrak s}=\big\{X \in {\mathfrak{sl}}(n, {\mathbb R})\colon X^{t}=X \big\}.
\end{gather*}

Now consider ${\mathfrak a}$ a maximal abelian subalgebra of ${\mathfrak s}$ given by all diagonal matrices  with va\-nishing trace. Also with canonical choices we have the following Iwasawa decomposition
\begin{gather*}
  {\mathfrak{sl}}(n, {\mathbb R})= {\mathfrak{so}}(n,{\mathbb R}) \oplus {\mathfrak a} \oplus {\mathfrak n}^{+} ,
\end{gather*}
where ${\mathfrak n}^{+}$ is the subalgebra of all upper triangular matrices with null entries on the diagonal.

Hence, the Iwasawa decomposition of the Lie group ${\rm SL}(n,\mathbb{R})$ is given by ${\rm SL}(n,\mathbb{R})= {\rm SO}(n,\mathbb{R})  \cdot A \cdot N^{+}$, where
\begin{gather*}
  A = \left\{ \textrm{diag}(a_1, \ldots, a_{n-1},a_n)
\colon a_1\cdot a_2\cdots   a_{n-1}\cdot a_n = 1
  \right\}
\end{gather*}
and
\begin{gather*}
  N^{+} = \left\{(y_{ij})\colon  y_{ii}=1, \, y_{ij} = 0,\, i>j, \, y_{ij} \in {\mathbb R}, \, 1<i<j\leq n
  \right\} .
\end{gather*}

Take the coordinate systems for $A$ and $N^{+}$, respectively, as
\begin{gather}\label{coordi}
  \phi(g) = (a_1,a_2,\ldots,a_{n-1}) \qquad \textrm{and} \qquad \psi(g) = (y_{12},y_{13}, \ldots,y_{n-1n}).
\end{gather}

Consider a smooth map $F\colon M \rightarrow {\rm SL}(n,\mathbb{R})/{\rm SO}(n,\mathbb{R})$. Then
\begin{gather*}
  F(x) = \big(F^1(x),F^2(x)\big) \in  A\times N^{+},\qquad \forall\, x \in M.
\end{gather*}
If $ v \in T_xM$ we have that $ \mu F_{x*}(v)$ is given by $\mu_{F}(x)(v)= \mu_{F^1}(x)(v) + \mu_{F^2}(x)(v) \in \mathfrak{a} \oplus \mathfrak{n}^{+}$. It is not dif\/f\/icult to see that
\begin{gather*}
  d^*\mu_{F}(x)(v)  = d^*\mu_{F^1}(x)(v) + d^*\mu_{F^2}(x)(v).
\end{gather*}

\subsection{Harmonic maps on canonical symmetric space}
In case of canonical Riemannian connection, $\beta\colon  \m \times \m \rightarrow \m$ is given by
\begin{gather*}
  \beta(X,Y) = \frac{1}{2}[X,Y]_m ,\qquad \textrm{for} \  X,Y \in \m.
\end{gather*}
Then from Example~\ref{iwasa-decomp} we have that $F$ is a harmonic map if and only if
\begin{gather*}
  d^*\mu_{F^1}(x) = 0\qquad  \textrm{and}\qquad d^*\mu_{F^2}(x)=0.
\end{gather*}
Writing $F_1$ and $F_2$ in the coordinates~(\ref{coordi}),
\begin{gather*}
  F^1(x)   =   \big(F^1_1(x),F^1_2(x),\ldots,F^1_{n-1}(x)\big) \qquad \textrm{and} \qquad
  F^2(x)   =   \big(F^2_{12}(x),F^2_{13}(x), \ldots,F^2_{n-1n}(x)\big).
\end{gather*}
we conclude that
\begin{Proposition}
  Let $F$ be a smooth map from $(M,g)$ into symmetric space ${\rm SL}(n,\mathbb{R})/{\rm SO}(n,\mathbb{R})$. Then~$F$ is harmonic if and only if the coordinates functions
  \begin{gather*}
   \big(F^1_1(x),F^1_2(x),\ldots,F^1_{n-1}(x)\big) \qquad \textrm{and} \qquad \big(F^2_{12}(x),F^2_{13}(x), \ldots,F^2_{n-1n}(x)\big)
  \end{gather*}
  are harmonic functions.
\end{Proposition}

As an immediate consequence we obtain the geodesics.
\begin{Corollary}
  Let $\gamma$ be a smooth curve in ${\rm SL}(n,\mathbb{R})/{\rm SO}(n,\mathbb{R})$ such that $\gamma(0) = \operatorname{Id}$ and
  \begin{gather*}
   \dot{\gamma}^1(0) ={\rm diag}\{a_1, \ldots, a_n\} \in \mathfrak{a} \qquad \mbox{and} \qquad \dot{\gamma}^2(0) = (n_{ij}) \in \mathfrak{n}.
  \end{gather*}
  Then $\gamma$ is geodesic if and only if
  \begin{gather*}
   \gamma^1(t) ={\rm diag}\{a_1t +1, \ldots, a_nt +1\} \qquad \mbox{and} \qquad \gamma^2(t)= (n_{ij}t).
  \end{gather*}
\end{Corollary}

\subsection[Harmonic maps with respect to $\alpha$-connections]{Harmonic maps with respect to $\boldsymbol{\alpha}$-connections}

In this case, the bilinear form $\beta\colon \m \times \m \rightarrow \m$ is given by
\begin{gather}\label{fisherconnection}
  \beta(X,Y) = \left(\dfrac{XY+YX}{2} - \dfrac{\operatorname{tr}(XY)}{n}\operatorname{Id}\right).
\end{gather}
First we need to compute  $\operatorname{tr}(\mu_F^*\nabla^G/H)$ in Theorem~\ref{mainteo}. In fact, for an orthonormal frame $\{e_1,\ldots, e_n\}$ on~$M$  we have that
\begin{gather*}
  \operatorname{tr}\big(\mu_F^*\nabla^G/H\big)
    =   \sum_{i=1}^{n}\beta(\mu_F(e_i),\mu_F(e_i))
    =   \sum_{i=1}^{n}\big(\mu_{F^1}(e_i)\cdot \mu_{F^1}(e_i)  - {\rm tr}(\mu_{F^1}(e_i)\cdot \mu_{F^1}(e_i))\operatorname{Id}\\
  \hphantom{\operatorname{tr}\big(\mu_F^*\nabla^G/H\big)}{}
    +   \mu_{F^1}(e_i)\cdot\mu_{F^2}(e_i)+ \mu_{F^2}(e_i)\cdot\mu_{F^1}(e_i)+ \mu_{F^2}(e_i)\cdot\mu_{F^2}(e_i)\big),
\end{gather*}
because the trace of any matrix in $\mathfrak{a}\oplus
\mathfrak{n}^{+}$ is null. Considering $\mu_{F^1}(e_i) \in
\mathfrak{a}$ and $ \mu_{F^2}(e_i) \in \mathfrak{n}^{+}$, Theorem
\ref{mainteo} assures that a smooth map $F\colon M \rightarrow {\rm
SL}(n,{\mathbb R})/{\rm SO}(n,{\mathbb R})$ is harmonic if and only if
\begin{gather}
  d^*\mu_{F^1} - \sum_{i=1}^{n} \left[\mu_{F^1}^2(e_i) -\frac{1}{n}\operatorname{tr}\big(\mu_{F^1}^2(e_i)\big)\operatorname{Id}\right]   =   0 \label{equation1}
\end{gather}
and
\begin{gather}
  d^*\mu_{F^2} - \sum_{i=1}^{n} \big[\mu_{F^1}(e_i)\cdot\mu_{F^2}(e_i)+ \mu_{F^2}(e_i)\cdot\mu_{F^1}(e_i)+ \mu_{F^2}^2(e_i)\big]  =  0. \label{equation2}
\end{gather}
As in dimension $2$ we have a better description, then we consider
the case of ${\rm SL}(2,{\mathbb R})/{\rm SO}(2,{\mathbb R})$.

Let $F\colon M^2 \rightarrow {\rm SL}(2,{\mathbb R})/{\rm SO}(2,{\mathbb R})$ be a smooth map. In coordinates,   $F = (F^1,F^2)=(F^1_1,F^1_{12})$. Thus $F$ is a harmonic map if it satisf\/ies the equations~(\ref{equation1}) and~(\ref{equation2}), namely, in a matrix context,
\begin{gather*}
  \left(
  \begin{matrix}
    \Delta F^1_{1} & 0 \\
    0 & -\Delta F^1_{1}
  \end{matrix}
  \right)
  -
  \left(
  \begin{matrix}
    \partial_x F^1_1 \partial_x F^1_{1} & 0\\
    0 & \partial_x F^1_1 \partial_x F^1_{1}
  \end{matrix}
  \right)
  -
  \left(
  \begin{matrix}
    \partial_yF^1_1\partial_yF^1_1 & 0 \\
    0 & \partial_yF^1_1\partial_y^1F_1
  \end{matrix}
  \right)
    \\
  \hphantom{  \left(
  \begin{matrix}
    \Delta F^1_{1} & 0 \\
    0 & -\Delta F^1_{1}
  \end{matrix}
  \right)}{}
  + \frac{1}{2}\left(
  \begin{matrix}
    \partial_x F^1_1 \partial_x F^1_{1} + \partial_yF^1_1\partial_yF^1_{1} & 0 \\
    0 & \partial_x F^1_1 \partial_x F^1_{1}+ \partial_yF^1_1\partial_yF^1_{1}
  \end{matrix}
  \right) =   0
\end{gather*}
and
\begin{gather*}
  \left(
  \begin{matrix}
    0 & \Delta F^2_{12}\\
    0 & 0
  \end{matrix}
  \right)
  -
  \left(
  \begin{matrix}
    0 & \partial_x F^1_1 \partial_x F^2_{12} \\
    0 & 0
  \end{matrix}
  \right)
  -
  \left(
  \begin{matrix}
    0 & -\partial_xF^1_1\partial_xF^2_{12}  \\
    0 & 0
  \end{matrix}
  \right)
 \\
\hphantom{  \left(
  \begin{matrix}
    0 & \Delta F^2_{12}\\
    0 & 0
  \end{matrix}
  \right)}{}
  -
  \left(
  \begin{matrix}
    0 & \partial_yF^1_1\partial_yF^2_{12}\\
    0 & 0
  \end{matrix}
  \right)
  -
  \left(
  \begin{matrix}
    0 & -\partial_yF^1_1\partial_yF^2_{12}\\
    0 & 0
  \end{matrix}
  \right)
    =   0.
\end{gather*}

A direct account shows that $\Delta F^1_{1} =0 $ and $\Delta F^1_{12} = 0$. In consequence,

\begin{Proposition}
  A smooth map $F\colon M^2 \rightarrow {\rm SL}(2,{\mathbb R})/{\rm SO}(2,{\mathbb R})$ is harmonic map if and only if coordinate functions $F^1_1$ and $F^2_{12}$ are real harmonic maps.
\end{Proposition}

\begin{Corollary}
  Let $\gamma\colon I \rightarrow {\rm SL}(2,{\mathbb R})/{\rm SO}(2,{\mathbb R})$ be a smooth curve such that $\gamma(0) = \operatorname{Id}$ and
  \begin{gather*}
    {\dot{\gamma}}^{1}(0) = {\rm diag}\{a,-a\}
  \end{gather*}
  and
  \begin{gather*}
    {\dot{\gamma}}^{2}(0) =
    \left(
    \begin{matrix}
      0 & n\\
      0 & 0
    \end{matrix}
    \right).
  \end{gather*}
  The curve $\gamma(t)=(\gamma^1_1(t),\gamma^2_{12}(t))$ is a geodesic if and only if $\gamma^1_1(t)= at+1$ and $\gamma^2_{12}(t) = nt$.
\end{Corollary}

The above results have not a good description for $n\geq3$, in fact we obtain systems of nonlinear partial and ordinary dif\/ferential equation.

\subsubsection[Geodesic in $N^+$]{Geodesic in $\boldsymbol{N^+}$}

Now our purpose is to study geodesics in  $N^+$. In the other hand,
Fernandes--San Martin~\cite{FernandesM/San-MartinL2} studied the
geodesics with respect to $\alpha$-connections in $A$.

\begin{Proposition}\label{geodesicN}
Let $\gamma(t)$ be a geodesic in ${\rm SL}(n,{\mathbb R})/{\rm SO}(n,{\mathbb R})$ with respect to $\alpha$-connec\-tion~\eqref{fisherconnection} that satisfies $\gamma(0) = \operatorname{Id}$,
$\gamma(t) \in N^+$ and $\dot{\gamma}(0) \in \mathfrak{n}$ given by
  \begin{gather*}
    \dot{\gamma}(0) =
    \left(
    \begin{matrix}
      0      & b_{12} & b_{13} & b_{14} &\ldots & b_{1n} \\
      0      &   0    & b_{23} & b_{24} &\ldots & b_{2n} \\
      0      &   0    &   0    & b_{34} &\ldots & b_{3n} \\
      \vdots & \vdots & \vdots & \ddots &\vdots &  0   \\
      0      &   0    &   0    & 0      &\ldots & b_{n-1n} \\
      0      &   0    &   0    & 0      &\ldots & 0
    \end{matrix}
    \right).
  \end{gather*}
  Then
  \begin{gather*}
    b_{ii+1}(t)   =   f(1)b_{ii+1}t,\qquad
    b_{ii+2}(t)   =   f(2)b_{ii+1}b_{i+1i+2}t^2,\\
    b_{ii+j}(t)   =   f(j)b_{ii+1}b_{i+1i+2}b_{i+2i+3}\cdots b_{(i+j-2)(i+j-1)}t^{j-1}, \qquad j>2,
  \end{gather*}
  where the function $f\colon \mathbb{N} \rightarrow \mathbb{N}$ is defined by recurrence as
  \begin{gather*}
    f(1)   =   1, \qquad
    f(2)   =   1, \qquad
    f(n)   =   \sum_{i=2}^{n} f(i-1)f((n+1)-i).
  \end{gather*}
\end{Proposition}

\begin{proof}
  First, it is clear from equation (\ref{equation2}) that
  \begin{gather*}
    b_{ii+1}(t)   =   f(1)b_{ii+1}t, \qquad
    b_{ii+2}(t)   =   f(2)b_{ii+1}b_{i+1i+2}t^2.
  \end{gather*}
  For $n\geq 3$ we will proof by induction. Suppose that we know the solution of all nilpotent matrix of dimension $n$. Taking a nilpotent matrix of dimension  $n+1$ we see that there are two nilpotent matrix of dimension $n$. Namely, the matrix $(b_{ij})_{1\leq i< j\leq n}$ and $(b_{ij})_{2\leq i< j\leq n+1}$, which we know the solutions. Then we need to f\/ind the solution of function $b_{1(n+1)}$. In fact, from equation~(\ref{equation2}) we see that this function satisf\/ies the  dif\/ferential equation
  \begin{gather*}
    \ddot{b}_{1(n+1)}(t) + \sum_{i=2}^{n} \dot{b}_{1i}(t)\dot{b}_{i(n+1)}(t) = 0.
  \end{gather*}
  We now know that
  \begin{gather*}
    \dot{b}_{1i}(t) = f(i-1)b_{12}b_{23} \cdots b_{(i-1)i}t^{(i-1)-1}
  \end{gather*}
  and
  \begin{gather*}
    \dot{b}_{i(n+1)}(t) = f((n+1)-i)b_{i(i+1)}b_{(i+1)(i+2)} \cdots b_{n(n+1)}t^{((n+1)-i)-1}.
  \end{gather*}
  Thus
  \begin{gather*}
    \ddot{b}_{1(n+1)}(t) + \sum_{i=2}^{n}[f(i-1)f((n+1)-i)]b_{12}b_{23}  \cdots b_{n(n+1)} t^{n-2}   =   0.
  \end{gather*}
  Therefore
  \begin{gather*}
    \ddot{b}_{1(n+1)}(t) + f(n)b_{12}b_{23}  \cdots b_{n(n+1)}t^{((n+1)i)-1} t^{n-2}  =  0.
  \end{gather*}
  Taking the initial condition we conclude that
  \begin{gather*}
    b_{1(n+1)}(t)=f(n)b_{12}b_{23} \cdots b_{n(n+1)}t^{((n+1)i)-1} t^{n}.\tag*{\qed}
  \end{gather*}
  \renewcommand{\qed}{}
\end{proof}

\begin{Remark}
  For example, if $\gamma(t)$ is a geodesic in ${\rm SL}(4,{\mathbb R})/{\rm SO}(4,{\mathbb R})$ with respect to $\alpha$-con\-nec\-tion~(\ref{fisherconnection}) such that $\gamma(t) \in N^+$, $\gamma(0) = \operatorname{Id}$ and
  \begin{gather*}
    \gamma(0) =
    \left(
    \begin{matrix}
      0      & b_{12} & b_{13} & b_{14} \\
      0      &   0    & b_{23} & b_{24}\\
      0      &   0    &   0    & b_{34}\\
      0      &   0    &   0    & 0
    \end{matrix}
    \right),
  \end{gather*}
  then Proposition \ref{geodesicN} give us that
  \begin{gather*}
    \gamma(t) =
    \left(
    \begin{matrix}
      0      & b_{12}t & b_{12}b_{23}t & b_{12}b_{23}b_{34}t^2 \\
      0      &   0    & b_{23}t & b_{23}b_{34}t\\
      0      &   0    &   0    & b_{34}t\\
      0      &   0    &   0    & 0
    \end{matrix}
    \right).
  \end{gather*}
\end{Remark}

\vspace{-3.5mm}

\subsection*{Acknowledgements}

\vspace{-1mm}

The authors gratefully acknowledge the many helpful suggestions of
the referees.
This work was partially supported by CNPq/Universal grant
n$^{\circ}$~$476024/2012{-}9$. The f\/irst author was partially
supported by Funda\c{c}\~{a}o Arauc\'{a}ria grant  n$^{\circ}$~$20134003$.

\vspace{-2.5mm}

\pdfbookmark[1]{References}{ref}
\LastPageEnding

\end{document}